\documentclass{amsart}
\usepackage{amssymb}
\usepackage{amsmath}
\usepackage{amscd}

\theoremstyle{plain}
\newtheorem{theorem}{Theorem}[section]
\newtheorem{lemma}[theorem]{Lemma}

\newtheorem{corollary}[theorem]{Corollary}
\newtheorem{proposition}[theorem]{Proposition}

\theoremstyle{definition}
\newtheorem{definition}[theorem]{Definition}

\newtheorem{examples}[theorem]{Examples}
\newtheorem{notation}[theorem]{Notation}
\newtheorem*{standing assumption}{Standing Assumption}
\newtheorem*{question}{Question}

\theoremstyle{remark}
\newtheorem{remark}[theorem]{Remark}
\newtheorem{remarks}[theorem]{Remarks}
\newtheorem*{ack}{Acknowledgment}

\begin{document}
\title[$C^*$-algebras from solenoids]
{$K$-theory of $C^*$-algebras from\\
one-dimensional generalized solenoids}

\author{Yi, Inhyeop}

\address{Department of Mathematics, 
University of Maryland, 
College park, MD, 20742}

\email{inhyeop{@}math.umd.edu}



\begin{abstract}
We compute the $K$-groups of $C^*$-algebras arising from
one-dimensional generalized solenoids.
The results show that Ruelle algebras from
one-dimensional generalized solenoids
are one-dimensional generalizations of
Cuntz-Krieger algebras. 
\end{abstract}

\maketitle

\section{Introduction}\label{s1} 
Ian Putnam and David Ruelle have developed a theory of
$C^*$-algebras for certain hyperbolic dynamical systems
(\cite{pu1, pu2, pusp, ru}).
These systems include
Anosov diffeomorphisms, topological Markov
chains and some examples of substitution tiling systems.
The corresponding $C^*$-algebras are modeled as
reduced groupoid $C^*$-algebras
for various equivalence relations.

This paper is concerned with $C^*$-algebras of
an orientable one-dimensional generalized solenoid
$(\overline{X},\overline{f})$,
where 
$\overline{X}$ has local canonical coordinates
which are contracting and expanding directions
for $\overline{f}$.
Na{\"{\i}}vely speaking,
Williams's orientable generalized solenoids
are higher dimensional analogues of
topological Markov chains (\cite{w1,w2}).
We consider the principal groupoids
of stable and unstable equivalence
on $(\overline{X},\overline{f})$,
denoted $G_s(\overline{X},\overline{f})$
and $G_u(\overline{X},\overline{f})$, respectively.
We give them topologies and Haar systems (\cite{pu1,pu2})
so that we may build
their reduced groupoid $C^*$-algebras
$S(\overline{X},\overline{f})$ and
$U(\overline{X},\overline{f})$, respectively,
as in \cite{re}.
The homeomorphism 
$\overline{f}\colon \overline{X}\to \overline{X}$
induces automorphism of $G_s(\overline{X},\overline{f})$
and $G_u(\overline{X},\overline{f})$,
and we form semi-direct products
$G_s\rtimes\mathbb{Z}$ and $G_u\rtimes \mathbb{Z}$.
Their groupoid $C^*$-algebras are denoted
$R_s(\overline{X},\overline{f})$
and $R_u(\overline{X},\overline{f})$, respectively,
and are called the {\it Ruelle algebras}
(\cite{pu2,pusp}).
In the case of topological Markov chains,
the Ruelle algebras are the Cuntz-Krieger algebras,
and the stable and unstable equivalence algebras are
the corresponding $AF$-subalgebras
of the Cuntz-Krieger algebras.

An important tool in the study of $C^*$-algebras
is $K$-theory.
Giordano, Herman, Putnam and Skau
showed that almost complete information about the orbit
structure of Cantor systems
is encoded by the $K$-theory
of their associated $C^*$-algebras
(\cite{gps,hps}). 
And Kirchberg and Phillips showed in their
recent papers (\cite{ki, ph}) that
nuclear, purely infinite, separable,
simple $C^*$-algebras are classified by
their $K$-theory.

In this paper, we compute
the $K$-groups of the unstable equivalence algebras
and the Ruelle algebras of $1$-solenoids
to answer the questions posed in \cite[\S4]{pu2}.
We show that the unstable equivalence algebra
of a $1$-solenoid $(\overline{X},\overline{f})$
with an adjacency matrix $M$
is strongly Morita equivalent to
the crossed product of a natural Cantor system
of $(\overline{X},\overline{f})$ by $\mathbb{Z}$
so that its $K_0$-group is order isomorphic to the dimension 
group of $M$ and its $K_1$-group is $\mathbb{Z}$.
Then we use the Pimsner-Voiculescu exact sequence,
the 
Universal Coefficient Theorem and 
Spanier-Whitehead duality to obtain
that the $K_0$-groups of Ruelle algebras
are isomorphic to
$\mathbb{Z}\oplus \{\Delta_M /\mathrm{Im}{}
(Id-\delta_M)\}$
and the $K_1$-groups are $\mathbb{Z}\oplus
\mathrm{Ker}{}(Id-\delta_M)$.
Thus $C^*$-algebras from
one-dimensional generalized solenoids are
one-dimensional analogues of the Cuntz-Krieger algebras.

The outline of the paper is as follow:
In section 2, we recall the axioms of
one-dimensional generalized solenoids and
their ordered group invariants.
In section 3, we review the definitions of
Smale spaces, and show that 
orientable one-dimensional solenoids are
Smale spaces.
Then we observe that
the stable equivalence algebras are
strongly Morita equivalent to inductive limit systems 
of $C^*$-algebras, and that the $K$-theory of
the unstable equivalence algebras
are determined by the adjacency matrices of
one-dimensional generalized solenoids.
In section 4,
we compute $K$-groups of unstable and stable
Ruelle algebras, and show that
they are $*$-isomorphic to each other by
the classification theorem of Kirchberg-Phillips.

\section{One-dimensional solenoids}
We review the properties of one-dimensional
generalized solenoids of Williams
which will be used in later sections.
As general references for the notions of
one-dimensional generalized solenoids
and their ordered group invariants
we refer to \cite{w1, w2, yi, yi2}.

\subsection*{One-dimensional generalized solenoids}
Let $X$ be a finite directed graph with vertex set
$\mathcal{V}$ and edge set $\mathcal{E}$,
and $f\colon X\to X$ a continuous map.
We define some axioms which
might be satisfied by $(X,f)$ (\cite{yi}).
\begin{enumerate}
\item[Axiom 0.]({\it Indecomposability})
$(X,f)$ is indecomposable.
\item[Axiom 1.]({\it Nonwandering})
All points of $X$ are nonwandering under $f$.
\item[Axiom 2.] ({\it Flattening}) 
There is $k\ge 1$ such that  for all $x\in X$ 
there is an open neighborhood $U$ of $x$
such that $f^k(U)$ is homeomorphic to
$(-\epsilon, \epsilon )$.
\item[Axiom 3.]({\it Expansion}) 
There are a metric $d$ compatible with the topology
and positive constants $C$ and $\lambda$ with 
$\lambda >1$ such that for all $n>0$ and 
all points $x,y$ on a common edge of $X$, 
if $f^n$ maps the interval $[x,y]$ into an edge,
then $d(f^nx,f^ny)\geq C\lambda^n d(x,y)$.  
\item[Axiom 4.] ({\it Nonfolding})
$f^n|_{X-\mathcal{V}}$ is locally one-to-one
for every positive integer $n$.
\item[Axiom 5.]({\it Markov})
$f(\mathcal{V})\subseteq \mathcal{V}$.
\end{enumerate}

Let $\overline{X}$ be the inverse limit space
$$
\overline{X}
=X\overset{f}{\longleftarrow}
X\overset{f}{\longleftarrow}\cdots
=\bigl\{(x_0,x_1,x_2,\dots )\in
       \prod_0^\infty X \, |\, f(x_{n+1})=x_n \bigr\},
$$
and $\overline{f}\colon \overline{X}\to \overline{X}$
the induced homeomorphism defined by
$$
(x_0,x_1,x_2,\dots )\mapsto 
(f(x_0),f(x_1),f(x_2),\dots )=(f(x_0),x_0,x_1,\dots).
$$

\begin{remark}\label{r2.1}
Williams' construction (\cite[6.2]{w2}) gives 
a (unique) measure $\mu_0$ for which
there is a constant $\lambda >1$ such that
$\mu_0 (X)=1$ and
$\mu_0 ( f(I) )=\lambda\mu_0(I)$
for every small interval $I\subset X$.
Define  $d_0(x_0,y_0)$
to be the measure of the smallest interval
from $x_0$ to $y_0$ in $X$, and
$$
d(x,y)=\sum\limits_{i=0}^{\infty}
\lambda^{-i}d_0(x_i,y_i)
$$
for $x=(x_0,x_1,x_2,\dots)$
and $y=(y_0,y_1,y_2,\dots)$ in $\overline{X}$.
Then $(\overline{X},d)$ is a compact metric space.
\end{remark}

Let $Y$ be a topological space and 
$g\colon Y\to Y$ a homeomorphism.
We call $Y$ a {\bf 1-dimensional generalized solenoid}
or {\bf $1$-solenoid} and $g$ a {\bf solenoid map}
if there exist a directed graph $X$
and a continuous map $f\colon X\to X$
such that $(X,f)$ satisfies all six Axioms
and $(\overline{X},\overline{f})$
is topologically conjugate to $(Y,g)$.
We call a point $x\in X$ a {\it non-branch point}
if $x$ has an open neighborhood
which is homeomorphic to an open interval,
and {\bf branch point} otherwise.
An {\bf elementary presentation} $(X,f)$ of a $1$-solenoid 
is such that $X$ is a wedge of circles and 
$f$ leaves the unique branch point of $X$ fixed.

\begin{proposition}[{\cite[5.2]{w2}}]\label{wep}
For each $1$-solenoid $(\overline{X}, \overline{f})$, 
there exists  an integer $m$ such that
$(\overline{X},\overline{f^m})$
has an elementary presentation.
\end{proposition}

Suppose that
$(X,f)$ is a presentation of a $1$-solenoid.
Since the inverse limit spaces
of $(X,f)$ and $(X,f^n)$
are homeomorphic (\cite{fo})
for every positive integer $n$,
for the purpose of computing invariants
of the space $\overline{X}$
there is no loss of generality
in replacing $(X,f)$ with $(X,f^n)$
where $n=m\cdot k$ is a positive integer such that
$(\overline{X}, \overline{f^m})$
has an elementary presentation $(Y,g)$
and for every $y\in Y$ there is an open set $U_y$
such that $g^k(U_y)$ is an open interval.
Hence we can assume that every point $x\in X$ has 
a neighborhood $U_x$ such that
$f(U_x)$ is an interval.

Recall that a continuous map $\gamma\colon [0,1]\to G$,
a directed graph,
is {\it orientation preserving}
if $e^{-1}\circ \gamma \colon I\to [0,1]$ is increasing
for every interval $I\subset [0,1]$
such that $\gamma(I)$ is a subset of a directed edge $e$.
A continuous map $\phi\colon G_1\to G_2$
between two directed graphs
is {\it orientation preserving}
if, for every orientation preserving map
$p \colon [0,1]\to G_1$,
the map $\phi\circ p \colon [0,1]\to G_2$ is
orientation preserving (\cite{fo}). 
 
When we can give a direction to each edge of $X$
so that the connection map $f\colon X\to X$ is
orientation preserving, we call
$({X},{f})$ an {\bf orientable presentation}.
For a $1$-solenoid $Y$ with a solenoid map $g$,
if there exists an orientable presentation $(X,f)$
then $Y$ is called an {\bf orientable} $1$-solenoid.

\begin{standing assumption}
In this paper, we always assume that $(X,f)$
is an orientable elementary presentation
such that every point $x\in X$ has 
a neighborhood $U_x$ such that
$f(U_x)$ is an interval.
\end{standing assumption}

\begin{notation}\label{4.n}
Suppose that $(X,f)$ is a presentation of a $1$-solenoid,
and that $\mathcal{E}=\{e_1,\dots,e_n\}$ is the edge set
of the directed graph $X$.
For each edge $e_i\in \mathcal{E}$,
we can give $e_i$ the {\it partition}
$\{I_{i,j}\}$, $1\le j \le l(i)$, such that
\begin{enumerate}
\item[(1)]
the initial point of $I_{i,1}$ is
the initial point of $e_i$,
\item[(2)]
the terminal point of $I_{i,j}$ is the initial
point of $I_{i,j+1}$ for $1\le j < l(i)$,
\item[(3)]
the terminal point of $I_{i,l(i)}$ is
the terminal point of $e_i$,
\item[(4)]
$f|_{\text{Int}{}I_{i,j}}$ is injective, and
\item[(5)]
$f(I_{i,j})=e_{i,j}^{s(i,j)}$ where
$e_{i,j}\in \mathcal{E}$,
$s(i,j)=1$ if the direction of $f(I_{i,j})$
agree with that of $e_{i,j}$, and
$s(i,j)=-1$ if the direction of $f(I_{i,j})$
is reverse to that of $e_{i,j}$.
\end{enumerate}
\end{notation}

The {\it wrapping rule}
$\check{f}\colon \mathcal{E}\to \mathcal{E}^*$
associated with $f$ is given by
$$
\check{f}\colon e_i \mapsto
e_{i,1}^{s(i,1)}\cdots e_{i,l(i)}^{s(i,l(i))},
$$
and the {\it adjacency matrix} $M$ of
$(\mathcal{E},\check{f})$ is given by
$$
M(i,k)=\#\{I_{i,j}\mid f(I_{i,j})=e_k^{\pm 1}\}. 
$$

\begin{remark}[{\cite[6.2]{w2}}]\label{r2.11}
The measure $\mu_0$ in remark \ref{r2.1} is given 
as follows:
Suppose that
$\lambda$ is the Perron-Frobenius eigenvalue
of the adjacency matrix $M$
and that $\mathbf{v}=(v_1, \dots, v_n)$ is
the corresponding Perron eigenvector such that
$\sum_{i=1}^{n}v_i=1$.
For edges $e_i, e_j$ of $X$ and
an interval $I$ of $e_i$ such that
$f^n(I)=e_{j}$ and $f^n|_{\text{Int}{}I}$ is injective,
let
$$
\mu_0(e_i)=v_i \text{ and }
\mu_0(I)=\lambda^{-n}v_{j}.
$$
Then $\mu_0$ is extended to a regular
Borel measure on $X$ by the standard procedure.
\end{remark}

\begin{theorem}[{\cite{am,ma,yi3}}]\label{tma}
Suppose that
$(\overline{X},\overline{f})$ is a $1$-solenoid.
Then there exists a uniquely ergodic flow $\phi$
whose phase space is $\overline{X}$.
\end{theorem}

Suppose that
$(X,f)$ is a presentation of a $1$-solenoid
and that $\mu_0$ is the measure
given on $X$ as in remark \ref{r2.11}.
For a measurable set $I$ in $X$,
we let
$U_n(I)=\{(x_0,\dots, x_n,\dots)\in \overline{X}
\mid x_n\in I\}$,
and  define a measure $\mu$ on $\overline{X}$ by
$$\mu\left( U_n(I)\right)=\mu_0(I).$$
Then $\mu$ is extended to a regular Borel measure
on $\overline{X}$ in the standard way.
We call this measure {\bf Williams measure}
of the flow $\phi$ on $\overline{X}$.
It is not difficult to verify that
$\mu$ is the unique $\phi$-invariant measure
on $\overline{X}$.

A closed subset $K$ of a phase space $Y$
of a flow $\phi$ is called a {\bf cross section}
if the mapping $\phi\colon  K\times \mathbb{R}\to Y$ 
defined by $(p,t)\mapsto p\cdot t$ is
a local homeomorphism onto $Y$.
The {\bf return time map} $r_k\colon  K\to K$ of
a cross section $K$ is defined by $x\mapsto y=x\cdot t_x$
where $x\in K$ and $t_x$ is the smallest positive number
such that $x\cdot t_x=y\in K$.

\begin{theorem}[\cite{hps, yi2}]\label{2.10}
Suppose that $(\overline{X},\overline{f})$ is
a $1$-solenoid with the corresponding
adjacency matrix $M$,
and that $(K,r_K)$ is a cross section with
the return time map of $\overline{X}$.
Then 
\begin{itemize}
\item[(1)]
$K_1(C(K)\times_{r_K}\mathbb{Z})=\mathbb{Z}$,
\item[(2)]
$K_0(C(K)\times_{r_K}\mathbb{Z})$
is order isomorphic to $\Delta_M$, and 
\item[(3)]
$K_0(C(K)\times_{r_K}\mathbb{Z})$ has a unique state.
\end{itemize}
\end{theorem}

\section{Smale spaces and $C^*$-algebras from solenoids}

\subsection*{Smale spaces $(${\cite{pu1, ru}}$)$}
Suppose that $(Y,d)$ is a compact metric space
and $\varphi$ is a homeomorphism of $Y$.
Assume that we have constants
$$
\text{ } 0< \lambda_0 <1,\text{ }\epsilon_0 > 0 
$$
and a continuous map
$$
(x,y)\in \left\{
(x,y)\in Y \times Y \mid d(x,y)\le 2\epsilon_0\right\}
\mapsto [x,y] \in Y
$$
satisfying the following:
$$
[x,x]=x,\quad
\left[ [x,y],z \right]=[x,z],\quad
\left[x,[y,z]\right]=[x,z],\quad
\left[ \varphi(x),\varphi(y)\right]=\varphi\left([x,y]\right)
$$
for $x,y,z\in Y$ whenever both sides of
the equation are defined.
For every $0< \epsilon \le \epsilon_0$ let
\begin{align*}
V^s(x,\epsilon)&=\left\{ y\in Y\mid [x,y]=y
\text{ and }d(x,y)< \epsilon\right\}\\
V^u(x,\epsilon)&=\left\{ y\in Y\mid [y,x]=y
\text{ and }d(x,y)< \epsilon\right\}.
\end{align*}
We assume that
\begin{align*}
d\left( \varphi(y),\varphi(z)\right)&\le
\lambda_0 d(y,z) \quad y,z\in V^s(x,\epsilon),\\
d\left( \varphi^{-1}(y),\varphi^{-1}(z)\right)&\le
\lambda_0 d(y,z) \quad y,z\in V^u(x,\epsilon).
\end{align*}
Then $(Y,d,\varphi)$ is called a {\bf Smale space}.

\subsection*{Groupoids $(${\cite{pu2,re}}$)$}
We refer to the work of Renault (\cite{re})
for the detailed theory of topological groupoids
and their associated $C^*$-algebras.
We give two examples of groupoids.

\begin{examples}[{\cite[1.2]{ru}}]
(1) {\it Equivalence relations}.
Suppose that $R$ is an equivalence relation
on a set $S$.
We give $R$ the following groupoid structure:
\begin{align*}
(s_1,t_1)\cdot (s_2,t_2)&=(s_1,t_2)
\text{ if } t_1=s_2 \text{ and}\\
(s,t)^{-1}&=(t,s).
\end{align*}

\noindent
(2) {\it Flows}.
Suppose that $S$ is a zero dimensional space and
$r\colon S \to S$ is a homeomorphism.
We consider the space $S\times\mathbb{R}$
with the equivalence relation,
$(s,\tau +1)\sim (r(s),\tau)$.
Let $\Sigma=S\times\mathbb{R}/\sim$
be the quotient space and define a flow
$\phi\colon \Sigma\times\mathbb{R}\to \Sigma$ by
$\phi_t(s,\tau)=\left[(s, t+\tau)\right]$
Give the following groupoid structure on
$\Sigma\times_{\phi} \mathbb{R}$:
\begin{align*}
(\sigma_1,t_1)\cdot(\sigma_2,t_2)&=(\sigma_1,t_1+t_2)
\text{ if $\sigma_2=\phi_{t_1}(\sigma_1)$ and}\\
(\sigma,t)^{-1}&=(\phi_t(\sigma),{-t}).
\end{align*}
\end{examples}

For a Smale space $(Y,d,\varphi)$,
define
$$
G_{s,0}=\left\{
(x,y)\in Y\times Y\mid y\in V^s(x,\epsilon_0)\right\}
\quad
G_{u,0}=\left\{
(x,y)\in Y\times Y\mid y\in V^u(x,\epsilon_0)\right\}
$$
and let
$$
G_s=\bigcup_{n=0}^\infty
\left(\varphi\times\varphi\right)^{-n}\left(G_{s,0}\right)
\quad
G_u=\bigcup_{n=0}^\infty
\left(\varphi\times\varphi\right)^{n}\left(G_{u,0}\right).
$$
Then $G_s$ and $G_u$ are equivalence relations on $Y$,
called {\it stable} and {\it unstable} equivalence.
Each $\left(\varphi\times\varphi\right)^{-n}
\left(G_{s,0}\right)$,
$\left(\varphi\times\varphi\right)^{-n}
\left(G_{u,0}\right)$
is given the relative topology of $Y\times Y$,
and $G_s$ and $G_u$
are given the inductive limit topology.
Then $G_s$ and $G_u$ are
locally compact Hausdorff principal groupoids.
The Haar systems $\{\mu_s^x\mid x\in Y\}$
and $\{\mu_u^x\mid x\in Y\}$
for $G_s$ and $G_u$, respectively,
are described in \cite[3.c]{pu2}.
The groupoid $C^*$-algebras of $G_s$ and $G_u$ are
denoted $S(Y,\varphi)$ and  $U(Y,\varphi)$,
respectively.

The map $\varphi\times \varphi$ acts
as an automorphism of $G_s$ and $G_u$.
We form the semi-direct products
\begin{align*}
G_s\rtimes \mathbb{Z}
&=\{(x,n,y)\mid n\in \mathbb{Z}\text{ and }
(\overline{f}^n(x),y)\in G_s\}\\
G_u\rtimes \mathbb{Z}
&=\{(x,n,y)\mid n\in \mathbb{Z}\text{ and }
(\overline{f}^n(x),y)\in G_u\}
\end{align*}
with groupoid operations
\begin{align*}
(x,n,y)\cdot (u,m,v)&=(x,n+m,v)
\text{ if } y=u \text{ and }\\
(x,n,y)^{-1}&=(y,-n,x).
\end{align*}
The product topology of $G_*\times \mathbb{Z}$
is transfered to $G_*\rtimes \mathbb{Z}$
by the bijective map 
$\eta\colon (x,y,n)\mapsto (x,n,\varphi(y))$.
And a Haar system on $G_*\rtimes \mathbb{Z}$
is given by $\mu_*^x\circ \eta^{-1}$
where $\mu_*^x$ is the Haar system on $G_*$.
The groupoid $C^*$-algebras
$C^*(G_s\rtimes \mathbb{Z})$
and $C^*(G_u\rtimes \mathbb{Z})$
are denoted $R_s(Y,\varphi)$ and $R_u(Y,\varphi)$
and are called the {\it Ruelle algebras}.

\begin{theorem}[{\cite{kps,pu1,pu2}}]\label{3.2}
Suppose that $(Y,\varphi)$ is a topologically mixing
Smale space. Then
\begin{itemize}
\item[(1)]
$S(Y,\varphi)$ and  $U(Y,\varphi)$ are
amenable, nuclear, separable and simple $C^*$-algebras, and
\item[(2)]
$R_s(Y,\varphi)$ and $R_u(Y,\varphi)$ are
amenable, non-unital, nuclear, purely infinite,
separable, simple and stable $C^*$-algebras.
\end{itemize}
\end{theorem}

For general properties of these $C^*$-algebras,
we refer to \cite{pu1,pu2,pusp}.

Suppose that $\left(\overline{X},\overline{f}\right)$
is a $1$-solenoid
with the metric $d$ given in remark \ref{r2.1}.
Let $\lambda_0=\epsilon_0=\frac{1}{\lambda}$
and define $[\; , \;]\colon 
\overline{X}\times \overline{X}\to\overline{X}$
by $[x,y]\mapsto z$
where
$z_0=x_0$ and $z_n$ is the unique element
contained in the $\lambda_0^{n+1}$-neighborhood of
$y_n$ such that $f^n(z_n)=x_0$.
Then it is not difficult to show that
$\left(\overline{X},\overline{f},d\right)$
satisfies the above conditions.
Therefore we have the following:
\begin{proposition}\label{p1}
One-dimensional generalized solenoids
are Smale spaces.
\end{proposition}


\subsection*{Stable equivalence algebras for 
$1$-solenoids}
Suppose that $G_s$ is
the stable equivalence groupoid
of a $1$-solenoid $(\overline{X},\overline{f})$
and that $S(\overline{X},\overline{f})$
is the corresponding groupoid algebra.
We first repeat the structural question
of Putnam (\cite[\S4]{pu2}).
For classical $1$-solenoids, we refer to \cite{d,pu1}.
\begin{question}
Can $S(\overline{X},\overline{f})$ be written as
an inductive limit?
\end{question}


\subsubsection*{Generalized transversals
$($\cite[\S3]{pusp}$)$}
Suppose that $G_s$ is the stable equivalence groupoid
of $(\overline{X},\overline{f})$,
that $U_p$ is the unstable equivalence class of
$p\in\overline{X}$ with the inductive limit topology
and that $g\colon U_p\to G_{s,0}$ is given by
$x\mapsto (x,x)$ for $x\in U_p$.
Let
$$
G_s(p)=\{(x,y)\in G_s \mid x,y\in U_p\}.
$$
A base for a topology on $G_s(p)$ is 
$$
\{U\cap s^{-1}\circ g(V^s)\cap r^{-1}\circ g(V^r)\mid
U\subset G_s, V^s,V^r\subset U_p
\text{ are open sets}\}.
$$
\begin{proposition}[{\cite[\S3]{pusp}}]
\begin{itemize}
\item[(1)]
$G_s(p)$ is an r-discrete, second countable,
locally compact, Hausdorff groupoid
with counting measure as Haar system.
\item[(2)]
$S(\overline{X},\overline{f})$
is strongly Morita equivalent to $C^*(G_s(p))$.
\end{itemize}
\end{proposition}

Now we choose $p$ to be
a fixed point of $\overline{f}$ such that
$\pi_k(p)$ is contained in the interior of an edge
$e\in \mathcal{E}$.
Since the orbits of $(\overline{X},\mathbb{R},\phi)$
are determined by the cofinality relation,
$x=(x_0,x_1,\dots)\in U_p$ if and only if
there is a positive integer $n=n(x)$ such that
$x_k \in e$ for every $k\ge n$.
Then
$(\overline{f}\times \overline{f})
\left(G_s(p)\right)=G_s(p)$.
Let
$$
G_{s,n}(p)=\{(x,y)\in G_{s,n}\mid x,y\in U_p\}
=\{(x,y)\in G_s(p)\mid f^n(x_0)=f^n(y_0)\}.
$$
Then $G_{s,n}(p)$ is
a compact open subset of $G_s(p)$,
and $G_{s,n}(p)^0 = G_s(p)^0=g(U_p)$.
Since $G_s(p)$ is $r$-discrete,
the range maps $r\colon G_s(p)\to G_s(p)^0$
and $r_n=r|_{G_{s,n}(p)}$ are local homeomorphisms.
Hence the Haar system of $G_s(p)$
restricted to $G_{s,n}(p)$
gives a Haar system for each $G_{s,n}(p)$.
Then we can express $C^*(G_s(p))$
as an inductive limit
$$
C^*(G_{s,1}(p))\to C^*(G_{s,2}(p))\to \cdots
\to C^*(G_{s,n}(p))\to \cdots.
$$

\subsection*{Unstable equivalence algebras}
Suppose that $\left(\overline{X},\overline{f}\right)$
is an orientable solenoid
and that $\phi$ is the flow on $\overline{X}$
given in theorem \ref{tma}.
Then there exists a cross section
with return time map $(K,r)$ such that
$\overline{X}$ is the suspension space of $(K,r)$.

\begin{lemma}[{\cite[p.59]{re}}]\label{l2.2}
The $C^*$-algebra of $(\overline{X},\mathbb{R},\phi)$
is isomorphic to $C(\overline{X})\times_{\phi}\mathbb{R}$.
\end{lemma}

\begin{proposition}[{\cite{mrw, pu2}}]
Suppose that $\left(\overline{X},\overline{f}\right)$
is an orientable solenoid,
and that $(Z,r)$ is
a cross section with the return time map
of the flow $\phi$.
Then
\begin{enumerate}
\item[(1)]
$U(\overline{X},\overline{f})
\simeq
C(\overline{X})\times_{\phi}\mathbb{R}$
and 
\item[(2)]
$C(\overline{X})\times_{\phi}\mathbb{R}$
is strongly Morita equivalent to
$C(K)\times_{r}\mathbb{Z}$.
\end{enumerate}
\end{proposition}
\begin{proof}
(1).
Suppose
$x=(x_0,x_1,\dots), y=(y_0,y_1,\dots) \in \overline{X}$
and $(x,y)\in G_u$.
Then $d\left( \overline{f}^n(x),\overline{f}^n(y)\right)
\to 0$ as $n\to -\infty$ implies
$d_0\left(x_n,y_n\right)\to 0$ as $n\to \infty$
and that there exists a $t\in \mathbb{R}$
such that $y=\phi_t(x)$.
Let
$\alpha\colon
(\overline{X},\mathbb{R},\phi) \to G_u$
be given by $(x,t) \mapsto \left(x,\phi_t(x)\right)$.
Then it is not difficult to see that
$\alpha$ is an isomorphism.
Therefore $U(\overline{X},\overline{f})$
is isomorphic to
$C(\overline{X})\times_{\phi}\mathbb{R}$
by lemma \ref{l2.2}.
And by the same argument
$C(K)\times_{r}\mathbb{Z}$ is isomorphic
to the groupoid $C^*$-algebra of $(K,\mathbb{Z},r)$.

\noindent
(2).
Since $\overline{X}$ is the suspension
of $(K,r)$, for every $x\in \overline{X}$
there exist unique $z_x\in K$ and $\tau_x\in [0,1)$
such that $x=\phi_{\tau_x}(z_x)$.
Define
$I=\left\{ (x,n-\tau_x)\mid
x\in \overline{X}, n\in \mathbb{Z}\right\}$,
and let $\mathcal{C}(I)$ be the completion of $C_c(I)$.
Then by the Theorem in \cite[\S4.a]{pu2}
$\mathcal{C}(I)$ is
a $C(\overline{X})\times_{\phi}\mathbb{R}$-$
C(K)\times_r\mathbb{Z}$ imprimitivity bimodule.
For completeness,
we write down the module structures and the inner
products.

\noindent{\it Module structures}.
Suppose that $\alpha\in C_c(I)$,
$g\in C_c(\overline{X},\mathbb{R},\phi)$ and
$h\in C_c(K,\mathbb{Z},r)$.
Then
\begin{align*}
\left(g\cdot \alpha\right) (x,n-\tau_x)
&=\int g(x,t)\cdot \alpha(\phi_t(x),n-\tau_x-t)
\,d\mu^{[x]}(t) \text{ and}\\
\left(\alpha\cdot h\right)(x,n-\tau_x)
&=\sum_m \alpha(x,m-\tau_x)\cdot h(r^m(z_x),n-m)
\end{align*}
give that 
$\mathcal{C}(I)$ is a left
$C(\overline{X})\times_{\phi}\mathbb{R}$
and right $C(K)\times_r\mathbb{Z}$ bimodule
with
$\left(\tilde{g}\cdot \tilde{\alpha}\right)\cdot\tilde{h}
=\tilde{g}\cdot (\tilde{\alpha}\cdot\tilde{h})$
for every
$\tilde{\alpha}\in \mathcal{C}(I)$,
$\tilde{g}\in C(\overline{X})\times_{\phi}\mathbb{R}$
and $\tilde{h}\in C(K)\times_r\mathbb{Z}$.

\noindent{\it Inner products}.
Define
$\langle \;,\;\rangle_{L}\colon
{C}_c(I)\times {C}_c(I) \to 
C_c(\overline{X},\mathbb{R},\phi)$
and
$\langle \;,\;\rangle_{R}\colon
{C}(I)\times{C}(I) \to C_c(K,\mathbb{Z},r)$ by
\begin{align*}
\langle \alpha,\beta\rangle_{L}(x,t)
&=\sum\alpha(x,m-\tau_x)\cdot 
\overline{\beta(x,m-\tau_x)} \text{ and}\\
\langle \alpha,\beta\rangle_{R}(z,k)
&=\int
\overline{\alpha\left(\phi_t(z),k-t\right)}
\cdot\beta\left(\phi_t(z),k-t\right)\,
d\mu^{[\phi_t(z)]}(t).
\end{align*}
\end{proof}

Then we have the following corollary from
propositions \ref{2.10}.

\begin{corollary}[{\cite{gps,yi2}}]
\begin{itemize}
\item[(1)]
$U(\overline{X},\overline{f})$
is a simple $C^*$-algebra.
\item[(2)]
$K_1\left(U(\overline{X},\overline{f})\right)
=\mathbb{Z}$.
\item[(3)]
$K_0\left(U(\overline{X},\overline{f})\right)$
is order isomorphic to $\Delta_M$
where $M$ is the adjacency matrix of
$(\overline{X},\overline{f})$.
\end{itemize}
\end{corollary}

Recall that the flow $\phi$ on $\overline{X}$ is
uniquely ergodic without rest point
(theorem \ref{tma}).
So $C(\overline{X})\times_{\phi}\mathbb{R}$
has the unique trace $\tau_{\mu}$ induced by
the Williams measure $\mu$ (\cite[3.3.10]{to}).
Thus $\tau_{\mu}^*$, the induced state
on $K_0(C(\overline{X})\times_{\phi}\mathbb{R})$,
is the unique state.


\begin{proposition}\label{3.6}
Suppose that $(\overline{X},\overline{f})$
is a $1$-solenoid and 
that $M$ is the corresponding adjacency matrix
with the normalized Perron eigenvector
$\mathbf{v}=(v_1,\dots,v_n)$.
Then
$$
\tau_{\mu}^*\left(
K_0(U(\overline{X},\overline{f}),
K_0(U(\overline{X},\overline{f}))_+
\right)=
\left< (\Delta_M,\Delta_M^+),\mathbf{v} \right>.
$$
\end{proposition}
\begin{proof}
Suppose that $\mathcal{E}_k=\mathcal{E}$ is
the edge set of the $k$th coordinate space of
$\overline{X}$.
Then by proposition \ref{2.10}
$$
\left(
K_0(U(\overline{X},\overline{f})),
K_0(U(\overline{X},\overline{f}))_+
\right)
\cong
\left(
\lim\limits_{\longrightarrow}
C(\mathcal{E}_k,\mathbb{Z}),
\lim\limits_{\longrightarrow}
C_+(\mathcal{E}_k,\mathbb{Z})
\right)
\cong
(\Delta_M,\Delta_M^+).
$$
For $g\in C(\mathcal{E}_k,\mathbb{Z})$,
$x=(x_0,\dots,x_k,\dots)\in \overline{X}$
with $x_k=e^{2\pi i s}\in e_i\in\mathcal{E}_k$
and the canonical projection
to the $k$th coordinate space
$\pi_k\colon \overline{X}\to X$,
define $g_k\in C(X_k,S^1)$
and $\tilde{g}\in C(\overline{X},S^1)$ by
$$
g_k\colon x_k\mapsto \exp(2\pi ig(e_i)s)
\text{ and }
\tilde{g}\colon x\to g_k\circ \pi_k (x).
$$
Then every $\tilde{g}$ is a unitary element
in $C(\overline{X})$,
and $K_0(U(\overline{X},\overline{f}))\cong
K_1(C(\overline{X}))$ is generated by $\tilde{g}$.
If we denote $g$ as $(g(e_1),\dots,g(e_n))$,
then by Theorem 2.2 of \cite{pa}
\begin{align*}
\tau_{\mu}^*(\tilde{g})
&=\frac{1}{2\pi i}\int_{\overline{X}}
\frac{\tilde{g}^{\prime}}{\tilde{g}}\,d\mu
=\int_{X_k} g^\prime \,d\mu_0
=\sum_{i=1}^n g(e_i)\mu_0(e_i)
=\sum_{i=1}^n g(e_i)v_i\\
&=\left< \left(g(e_1),\dots,g(e_n)\right),
\mathbf{v} \right>.
\end{align*}
\end{proof}

The above proposition refines Theorem 2.2 of \cite{pa}
that
$$\tau_{\mu}^*(K_0(C(\overline{X})\times_{\phi}\mathbb{R}))
=\left< A_{\mu},\check{H}^1(\overline{X})\right>.
$$

\begin{corollary}[{\cite{bl}}]
If $p$ and $q$ are projections in 
$M_{\infty}\left( C(\overline{X})
\times_{\phi}\mathbb{R}\right)$ such that
$\tau_{\mu}(p)< \tau_{\mu}(q)$, then 
$p$ is equivalent to a subprojection of $q$.
\end{corollary}


\begin{lemma}[{\cite{pu0}}]\label{lsr}
$C(K)\times_r\mathbb{Z}$ has
real rank zero and topological stable rank one.
\end{lemma}

Since  $C(\overline{X})\times_{\phi}\mathbb{R}$
and $C(K)\times_r\mathbb{Z}$ are separable algebras,
they have strictly positive elements.
So strong Morita equivalence of
$C(\overline{X})\times_{\phi}\mathbb{R}$
and $C(K)\times_r\mathbb{Z}$
implies that they are stably isomorphic, i.e.,
$\{C(\overline{X})\times_{\phi}\mathbb{R}\}\otimes \mathcal{K}$
is $*$-isomorphic to
$\{C(K)\times_r\mathbb{Z}\}\otimes \mathcal{K}$
where $\mathcal{K}$ is the algebra of compact operators
on a separable Hilbert space.
Therefore we have the following proposition.

\begin{proposition}
$U(\overline{X},\overline{f})$
has real rank zero and topological stable rank one.
\end{proposition}

\section{Ruelle algebras for solenoids}
We compute $K$-groups of Ruelle algebras
for $1$-solenoids
to show that they are $*$-isomorphic.
\subsection*{Unstable equivalence Ruelle algebras}
Suppose that $(\overline{X},\overline{f})$
is an oriented $1$-solenoid and that
$G_u\simeq (\overline{X},\mathbb{R},\phi)$
is the unstable equivalence groupoid on $\overline{X}$.
Recall that for $x,y\in \overline{X}$
such that $y=\phi_t(x)$, $t\in \mathbb{R}$,
we have
$\overline{f}^{-1}(y)
=\phi_{t\lambda^{-1}}\circ\overline{f}^{-1}(x)$.

\begin{definition}[{\cite[\S4]{pu2}}]
Let $\alpha_u$ be an automorphism on
$U(\overline{X},\overline{f})$ defined by
$$
\alpha_u(g)(x,t)
=\lambda^{-1}g(\overline{f}^{-1}(x),t\lambda^{-1})
\text{ for } g\in C_c(\overline{X},\mathbb{R},\phi)
\text{ and } (x,t)\in (\overline{X},\mathbb{R}).
$$
The {\it unstable equivalence Ruelle algebra}
$R_u(\overline{X},\overline{f})$ is the crossed product
$$
R_u(\overline{X},\overline{f})
=U(\overline{X},\overline{f})\times_{\alpha_u}\mathbb{Z}
=\left(C(\overline{X})\times_{\phi}\mathbb{R}\right)
\times_{\alpha_u}\mathbb{Z}.
$$
\end{definition}

\begin{remarks}
\begin{itemize}
\item[(1)]
Let $A$ be an $n\times n$ integer matrix
and $\Delta_A$ the dimension group of $A$.
The {\it dimension group automorphism}
$\delta_A$ of $A$ is the restriction of $A$ to $A$
so that $\delta_A(\mathbf{v})=A\mathbf{v}$
(\cite[7.5.1]{lm}). 
Then 
$\Delta_A/ \mathrm{Im}{}(Id-\delta_A)$
is isomorphic to $\mathbb{Z}^n / (Id - A)\mathbb{Z}^n$.
\item[(2)]
For $g\in C(\mathcal{E}_k,\mathbb{Z})$,
let 
$g_k\in C(X_k,S^1)$ be
as in the proof of proposition \ref{3.6}.
The wrapping rule
$\check{f}\colon \mathcal{E}_{k+1}\to \mathcal{E}_k$
induces a map
$f^*\colon  C(\mathcal{E}_k,\mathbb{Z})
\to C(\mathcal{E}_{k+1},\mathbb{Z})$
by $g\mapsto g\circ \check{f}$ where
$(g\circ \check{f})(e)=\sum\limits_{i=1}^{j}g(e_i)$
such that $\check{f}(e)=e_1\cdots e_j$.
Then 
$
g_k\circ f\circ \pi_k$ is homotopic to
$(g\circ {f}^*)_{k+1} \circ \pi_{k+1}$
(\cite[3.6]{yi2}).
\end{itemize}
\end{remarks}

\begin{proposition}\label{3.12}
Suppose that $(\overline{X},\overline{f})$
is a $1$-solenoid with
the adjacency matrix $M$
and corresponding dimension group automorphism
$\delta_M$.
Then
$$
K_0(R_u(\overline{X},\overline{f}))\cong 
\mathbb{Z}\oplus
\{\Delta_M/ \mathrm{Im}{}(Id - \delta_M )\}
\text{ and }
K_1(R_u(\overline{X},\overline{f}))
\cong\mathbb{Z} \oplus \mathrm{Ker}(Id-\delta_M ).
$$
\end{proposition}
\begin{proof}
We have the following Pimsner-Voiculescu exact sequence.
$$
\begin{CD}
K_0(U(\overline{X},\overline{f}))@>1-{\alpha_u}_*>>
K_0(U(\overline{X},\overline{f}))@>\iota_*>>
K_0(R_u(\overline{X},\overline{f}))\\
@AAA @. @VVV\\
K_1(R_u(\overline{X},\overline{f}))
@<<\iota_*< K_1(U(\overline{X},\overline{f}))
@<<1-{\alpha_u}_*< K_1(U(\overline{X},\overline{f}))
\end{CD}
$$
We consider
${\alpha_u}_*\colon
K_0(U(\overline{X},\overline{f}))
=K_0\left(
C(\overline{X})\times_{\phi}\mathbb{R}
\right)\to
K_0\left(
C(\overline{X})\times_{\phi}\mathbb{R}
\right)$
as the automorphism
$\hat{\alpha}_{u*}\colon
K_1(C(\overline{X}))\to K_1(C(\overline{X}))$
given by the Thom isomorphism of Connes.
Define 
$\beta \colon C(\overline{X})
\to C(\overline{X})$ by
$h \mapsto h\circ \overline{f}^{-1}$
for $h\in C(\overline{X})$.
Then the induced automorphism
$\beta_*\colon
K_1(C(\overline{X})) \to K_1(C(\overline{X}))$
is the required isomorphism.

For $g\in C(\mathcal{E}_k,\mathbb{Z})$,
let $\tilde{g}\in C(\overline{X},S^1)$
be the induced unitary element
as in the proof of proposition \ref{3.6}.
Then 
$\beta^{-1}(\tilde{g})=\tilde{g}\circ \overline{f}
=g_k\circ\pi_k\circ \overline{f}
=g_k\circ f\circ \pi_k$ is homotopic to
$(g\circ {f}^*)_{k+1} \circ \pi_{k+1}$.
Hence if we denote $g$ as
$(g(e_1),\dots,g(e_n))\in \mathbb{Z}^n$,
then $g\circ {f}^*$ is given by $M g$
and the induced automorphism
$\beta_*^{-1}\colon
K_1(C(\overline{X})) \to K_1(C(\overline{X}))$
is the dimension group automorphism $\delta_M$
of the adjacency matrix $M$.
Therefore 
$\beta_*$ is the inverse of $\delta_M$,
and $1-\alpha_{u*} \colon
K_0(U(\overline{X},\overline{f}))
\to K_0(U(\overline{X},\overline{f}))$
is the same as
$Id-\delta_M^{-1}\colon \Delta_M\to \Delta_M$.

Since $K_1(U(\overline{X},\overline{f}))$
is isomorphic to $\mathbb{Z}$,
$\alpha_{u*} \colon \mathbb{Z} \to \mathbb{Z}$
is trivially the identity map.
Thus the six-term exact sequence is divided into
the following two short exact sequences;
$$
0\to \Delta_M / \mathrm{Im}(Id-\delta_M^{-1})
\longrightarrow K_0(R_u(\overline{X},\overline{f}))
\longrightarrow \mathbb{Z}\to 0
$$
and
$$
0\to \mathbb{Z}\longrightarrow
K_1(R_u(\overline{X},\overline{f}))
\longrightarrow \mathrm{Ker}(Id-\delta_M^{-1})
\to 0.
$$
Therefore we conclude that
\begin{align*}
K_0(R_u(\overline{X},\overline{f}))
&\cong  \mathbb{Z}\oplus
\{\Delta_M/ \mathrm{Im}{}(Id - \delta_M^{-1} )\}
=\mathbb{Z}\oplus
\{\Delta_M/ \mathrm{Im}{}(Id - \delta_M)\}
\text{ and }\\
K_1(R_u(\overline{X},\overline{f}))
&\cong \mathbb{Z} \oplus \mathrm{Ker}(Id-\delta_M^{-1})
=\mathbb{Z} \oplus \mathrm{Ker}(Id-\delta_M).
\end{align*}
\end{proof}

\begin{examples}
(1).
Suppose that $X$ is the unit circle
and that $f\colon X\to X$ is given by
$z\mapsto z^n$, $n\ge 2$.
Then the adjacency matrix is $(n)$,
$K_0(U(\overline{X},\overline{f}))
=\mathbb{Z}[\frac{1}{n}]$
and $K_1(U(\overline{X},\overline{f}))
=\mathbb{Z}$.
Since $\delta_{(n)}^{-1}$ is multiplication by
$\frac{1}{n}$,
we have
$K_0(R_u(\overline{X},\overline{f}))
=\mathbb{Z}\oplus \mathbb{Z}_{n-1}$
and
$K_1(R_u(\overline{X},\overline{f}))
=\mathbb{Z}$.
See \cite{d, ls} for details.

(2).
Suppose that $Y$ is a wedge of two circles $a$ and $b$
and that $g\colon Y\to Y$ is given by
$a\mapsto aab$ and $b\mapsto ab$.
Then the adjacency matrix is
$M=\begin{pmatrix} 2&1\\1&1\end{pmatrix}$.
So $K_0(U(\overline{Y},\overline{g}))
=\mathbb{Z}\oplus \mathbb{Z}$
and $K_1(U(\overline{Y},\overline{g}))=\mathbb{Z}$.
Since $1-{\alpha_u}_*\colon \mathbb{Z}\oplus \mathbb{Z}
\to \mathbb{Z}\oplus \mathbb{Z}$
is an isomorphism, 
we 
obtain
$K_0(R_u(\overline{Y},\overline{g}))
=K_1(R_u(\overline{Y},\overline{g}))
=\mathbb{Z}$.
\end{examples}

\subsection*{Stable equivalence Ruelle algebras}
We use $K$-theoretic duality of the Ruelle algebras
and the Universal Coefficient Theorem 
to compute $K$-groups of
$R_s(\overline{X},\overline{f})$.

\begin{remark}[{\cite{rs}}]
Let $\mathcal{N}$ be the category of
separable nuclear $C^*$-algebras
which contains the separable Type I $C^*$-algebras
and is closed under strong Morita equivalence,
inductive limits, extensions,
and crossed products by $\mathbb{Z}$
and by $\mathbb{R}$.
Then it is not difficult to verify that
unstable and stable equivalence Ruelle algebras of
$1$-solenoids are contained in $\mathcal{N}$.
\end{remark}

\begin{proposition}[{\cite[5.c]{pu2}}]
Suppose that $(\overline{X},\overline{f})$
is a $1$-solenoid.
Then $R_s(\overline{X},\overline{f})$
is dual to $R_u(\overline{X},\overline{f})$
so that $K_*(R_s(\overline{X},\overline{f}))$
is isomorphic to
$K^{*+1}(R_u(\overline{X},\overline{f}))$.
\end{proposition}

\begin{proposition}[{\cite[1.19]{rs}}]
Suppose that
$(\overline{X},\overline{f})$ is a $1$-solenoid.
Then there are short exact sequences
\begin{gather*}
0\to 
\mathrm{Ext}_{\mathbb{Z}}^1
(K_0(R_u(\overline{X},\overline{f})),\mathbb{Z})
\to K^1(R_u(\overline{X},\overline{f}))\to 
\mathrm{Hom}
(K_1(R_u(\overline{X},\overline{f})),\mathbb{Z})
\to 0\\
0\to 
\mathrm{Ext}_{\mathbb{Z}}^1
(K_1(R_u(\overline{X},\overline{f})),\mathbb{Z})
\to K^0(R_u(\overline{X},\overline{f}))\to 
\mathrm{Hom}
(K_0(R_u(\overline{X},\overline{f})),\mathbb{Z})
\to 0
\end{gather*}
\end{proposition}

Hence $K$-groups of the stable equivalence Ruelle
algebra
are determined by
$\mathrm{Ext}$- 
and $\mathrm{Hom}$-groups of 
$K_*(R_u(\overline{X},\overline{f}))$.
%
%
Transform $Id - M$ to the Smith form
$$
\begin{pmatrix}
d_1& & & \\ &d_2& &\\ & &\ddots & \\ & & &d_n
\end{pmatrix}
$$
where $d_i\ge 0$ and $d_i$ divides $d_{i+1}$
(\cite[\S7.4]{lm}).
Then
$\Delta_M/\mathrm{Im}{}(Id-\delta_M)$
is isomorphic to
$\oplus_{i=1}^{n}\,\mathbb{Z}/{d_i \mathbb{Z}}$,
and the dimension of $\mathrm{Ker}(Id -\delta_M)$
is equal to the number of zeros in the diagonal 
of the Smith form.
Suppose $d_1=\dots=d_m=0$ and $d_{m+1}\ne 0$.
Then we have
\begin{align*}
\mathrm{Ext}_{\mathbb{Z}}^1
(K_0(R_u(\overline{X},\overline{f})),\mathbb{Z})
&=\mathrm{Ext}_{\mathbb{Z}}^1
(\mathbb{Z}_{d_{m+1}}\oplus\cdots
\oplus\mathbb{Z}_{d_{n}},\mathbb{Z})
=\mathbb{Z}_{d_{m+1}}\oplus\cdots
\oplus\mathbb{Z}_{d_{n}} \text{ and}\\
\mathrm{Hom}(K_1
(R_u(\overline{X},\overline{f})),\mathbb{Z})
&=\mathbb{Z}^{m+1}.
\end{align*}
Hence we have
\begin{align*}
K^1(R_u(\overline{X},\overline{f}))
&\cong
\mathrm{Hom}(K_1(R_u
(\overline{X},\overline{f})),\mathbb{Z}) \oplus
\mathrm{Ext}_{\mathbb{Z}}^1
(K_0(R_u(\overline{X},\overline{f})),\mathbb{Z})\\
&=\mathbb{Z}\oplus \mathbb{Z}^{m}\oplus
\mathbb{Z}_{d_{m+1}}\oplus\cdots\mathbb{Z}_{d_{n}}\\
&\cong\mathbb{Z}\oplus
\{\Delta_{M}/\mathrm{Im}{}(Id - \delta_{M})\}.
\end{align*}

Recall that
$K_1(R_u(\overline{X},\overline{f}))
=\mathbb{Z}\oplus \mathrm{Ker}(Id-\delta_M)$
is a torsion-free subgroup of $\mathbb{Z}^{n+1}$.
Thus we have
$\mathrm{Ext}_{\mathbb{Z}}^1
(K_1(R_u(\overline{X},\overline{f})),\mathbb{Z})
=0$ and 
$$
K^0(R_u(\overline{X},\overline{f}))
\cong\mathrm{Hom}
(K_0(R_u(\overline{X},\overline{f})),\mathbb{Z}).
$$
Then 
$K_0(R_u(\overline{X},\overline{f}))
\cong \mathbb{Z}
\oplus_{i=1}^{n}\mathbb{Z}/{d_i \mathbb{Z}}$
implies 
$$
\mathrm{Hom}
(K_0(R_u(\overline{X},\overline{f})),\mathbb{Z})
\cong
\mathrm{Hom}(\mathbb{Z}
\oplus_{i=1}^{n}\mathbb{Z}/{d_i \mathbb{Z}},\mathbb{Z})
\cong \mathbb{Z}\oplus_{i=1}^{m}\mathbb{Z}
\cong \mathbb{Z}\oplus\mathrm{Ker}(Id-\delta_{M}).
$$
Therefore we conclude that:
\begin{proposition}\label{3.16}
Suppose that
$(\overline{X},\overline{f})$ is a $1$-solenoid.
Then
$$
K_0(R_s(\overline{X},\overline{f}))
\cong \mathbb{Z}\oplus\{\Delta_{M} /
\mathrm{Im}{}(Id - \delta_{M})\}
\text{ and }
K_1(R_s(\overline{X},\overline{f}))\cong
\mathbb{Z}\oplus\mathrm{Ker}(Id-\delta_{M}).
$$
\end{proposition}

\begin{remark}
The isomorphisms in proposition \ref{3.16}
are {\it  unnatural} as the short exact sequences
in the Universal Coefficient Theorem
split unnaturally.
\end{remark}

Recall that the unstable and stable equivalence
Ruelle algebras of a $1$-solenoid are
nuclear, purely infinite, separable,
simple and stable $C^*$-algebras (proposition \ref{3.2}).
Then the classification theorem of Kirchberg-Phillips
implies the following proposition.

\begin{proposition}
$R_u(\overline{X},\overline{f})$
is $*$-isomorphic to
$R_s(\overline{X},\overline{f})$.
\end{proposition}

\begin{ack}
I express my deep gratitude to
Dr.\;M.\;Boyle and Dr.\;J.\;Rosenberg at UMCP
and Dr.\;I.\;Putnam at University of Victoria, Canada,
for their encouragement and useful discussions.
The $[\,,\,]$-function for $1$-solenoids was
suggested by Dr.\;Putnam.
By kind permission, I presented his definition.
\end{ack}

\bibliographystyle{amsplain}

\end{document}